\documentclass[10pt]{article}

\usepackage{amsmath}
\usepackage{amsfonts}
\usepackage{amssymb}
\usepackage{framed}
\usepackage{mathtools}
\usepackage{enumerate}
\usepackage{mathrsfs}
\usepackage[hidelinks]{hyperref}
\usepackage{enumitem}
\usepackage{mathrsfs}
\usepackage{scrextend}
\usepackage{amsthm}
\usepackage{bbm}
\usepackage{thmtools}
\usepackage{thm-restate}
\usepackage{mathabx}
\usepackage[margin=1.5in]{geometry}
\usepackage{rotating}

\usepackage{xcolor}
\hypersetup{
    colorlinks,
    linkcolor={red!50!black},
    citecolor={blue!50!black},
    urlcolor={blue!80!black}
}

\usepackage{tikz-cd}

\theoremstyle{plain}
\newtheorem{thm}{Theorem}[section]
\newtheorem{lem}[thm]{Lemma}
\newtheorem{prop}[thm]{Proposition}
\newtheorem{cor}[thm]{Corollary}

\theoremstyle{definition}
\newtheorem{defn}[thm]{Definition}

\newtheorem{conj}[thm]{Conjecture}

\theoremstyle{remark}
\newtheorem{rem}[thm]{Remark}

\newtheorem*{notn}{Notation}

\tikzset{
  symbol/.style={
    draw=none,
    every to/.append style={
      edge node={node [sloped, allow upside down, auto=false]{$#1$}}}
  }
}

\newcommand\restr[2]{{
	\left.\kern-\nulldelimiterspace
	#1
	\vphantom{\big|}
	\right|_{#2}
	}}

\newcommand{\ep}{\varepsilon}

\newcommand{\ch}[1]{\widecheck{{#1}}}

\newcommand{\NL}{\operatorname{NL}}
\newcommand{\codim}{\operatorname{codim}}

\newcommand{\GL}{\operatorname{GL}}

\newcommand{\SL}{\operatorname{SL}}

\newcommand{\Sp}{\operatorname{Sp}}
\newcommand{\HL}{\operatorname{HL}}

\usepackage{amsmath,calligra,mathrsfs}

\usepackage[cal=boondoxo]{mathalfa}

\DeclareMathOperator{\sheafhom}{\mathcal{H \kern -1pt o \kern -2pt m}}
\DeclareMathOperator{\sheafend}{\mathcal{E \kern -1pt n \kern -2pt d}}
\DeclareMathOperator{\sheafaut}{\mathcal{A \kern -1pt u \kern -2pt t}}

\title{On The Complexity of Atypical Special Points}
\author{David Urbanik}

\begin{document}

\maketitle

\begin{abstract}
Given an integral variation of Hodge structure $\mathbb{V}$ on a complex algebraic variety $S$, polarized by some bilinear form $Q : \mathbb{V} \otimes \mathbb{V} \to \mathbb{Z}$, it is believed that the set $\mathcal{A}^{\textrm{iso}}_{0} \subset S(\mathbb{C})$ of isolated atypical special points associated to $(\mathbb{V}, Q)$ forms a finite set. Here we show that the number of such points $s$ is $O(Q(t_{s}, t_{s})^{\ep})$ for any $\ep > 0$, where $t_{s}$ is a minimal integral Hodge tensor defining $s$ (in an appropriate sense). This resolves a conjecture of Grimm and Monnee.
\end{abstract}

\tableofcontents

\subsection{Conventions}

All Mumford-Tate groups in this paper are special Mumford-Tate groups. This means that if $h : \mathbb{S} \to \GL(V)_{\mathbb{R}}$ is a Hodge structure, with $\mathbb{S}$ the Deligne torus, its Mumford-Tate group is the $\mathbb{Q}$-Zariski closure of $h(\mathbb{U})$, where $\mathbb{U} \subset \mathbb{S}$ is the ``circle'' subtorus defined in \cite[\S I.A]{GGK}. 

\section{Introduction}

\subsection{General Background}

Given a polarized integral variation of Hodge structure $(\mathbb{V}, Q)$ on a smooth complex algebraic variety $S$, one obtains a set $\mathcal{S}$ of \emph{special subvarieties} of $S$. They are defined as follows. Given any point $s_{0} \in S(\mathbb{C})$, and any tensor $t \in \bigoplus_{a \geq 0} \mathbb{V}^{\otimes a}_{s_{0}}$, there is a locus
\[ \HL(S,t) = \left\{ s \in S(\mathbb{C}) : \textrm{some parallel translate of }t\textrm{ to } \bigoplus_{a \geq 0} \mathbb{V}^{\otimes a}_{s}\textrm{ is Hodge} \right\} . \]
A theorem of Cattani-Deligne-Kaplan \cite{CDK} guarantees that $\HL(S,t)$ is an algebraic subvariety. Then $\mathcal{S}$ is the set of subvarieties of $S$ that are obtained as irreducible components of $\HL(S,t)$ for some $t$. Note that $S \in \mathcal{S}$ always by taking $t = 0$. 

\begin{defn}
We say a special subvariety $Z \subset S$ is strict if it is strictly contained in $S$.
\end{defn}

The Zilber-Pink conjecture studies the subset $\mathcal{A} \subset \mathcal{S}$ of maximal atypical special subvarieties. To define the term \emph{atypical} one can use \emph{Hodge data}. We recall that, following \cite{atyp}, to any variation of Hodge structure $\mathbb{W}$ on an irreducible algebraic variety $Y$ one can associate a pair $(\mathbf{G}_{Y,\mathbb{W}}, D_{Y,\mathbb{W}})$ called the Hodge datum of $(Y, \mathbb{W})$. When $\mathbb{W}$ is understood we write simply $(\mathbf{G}_{Y}, D_{Y})$. In this case, we are interested in the situations $(Y, \mathbb{W}) = (Z, \restr{\mathbb{V}}{Z})$, with $Z \subset S$ a special subvariety, and their associated Hodge data $(\mathbf{G}_{Z}, D_{Z})$. When $Z = S$ we write $(\mathbf{G}, D) := (\mathbf{G}_{S}, D_{S})$. 
\begin{defn}
A special subvariety $Z \in \mathcal{S}$ is said to be atypical if $\codim_{S} Z < \codim_{D_{S}} D_{Z}$. 
\end{defn}
One then has, following \cite[Conj 2.5]{BKU}:
\begin{conj}[Zilber-Pink]
\label{ZPconj}
Given any polarizable variation of Hodge structure $\mathbb{V}$ on a complex algebraic variety $S$, there are only finitely many maximal atypical special subvarieties of $S$ for $\mathbb{V}$ under inclusion.
\end{conj}
To simplify matters, let us work in the additional setting where the adjoint group $\mathbf{G}^{\textrm{ad}}_{S}$ is a $\mathbb{Q}$-simple group; this assumption is usually satisfied in practice, for instance whenever the image of monodromy acting on the fibres of $\mathbb{V}$ is sufficiently large. In this case \emph{geometric Zilber-Pink} \cite[Thm 3.1]{BKU} gives the following result:
\begin{thm}[BKU]
\label{geoZP}
When $\mathbf{G}^{\textrm{ad}}_{S}$ is $\mathbb{Q}$-simple, the Zariski closure $Z_{\textrm{poscl}}$ in $S$ of the collection of all atypical special subvarieties of positive period dimension is contained in a finite union $Z_{1} \cup \cdots \cup Z_{k}$ of strict special subvarieties of $S$.
\end{thm}

\begin{proof}
Apply \cite[Thm. 3.1]{BKU} to each component $Z$ of $Z_{\textrm{poscl}}$. If condition (a) in loc. cit. occurs then we are done, otherwise (b) occurs and $Z \subsetneq S$ since $\mathbf{G}^{\textrm{ad}}_{S}$ is $\mathbb{Q}$-simple and (b) gives a non-trivial splitting of $\mathbf{G}_{Z}$. Then $Z$ is contained in a special subvariety since $\mathbf{G}_{Z} \subsetneq \mathbf{G}_{S}$.
\end{proof}

\begin{rem}
Note that the conclusion of \autoref{geoZP} does not guarantee that the $Z_{i}$ are atypical, which would follow from \autoref{ZPconj}.
\end{rem}

\noindent Let us explain the term \emph{positive period dimension}. Because the polarizing form induces a positive-definite bilinear form on the lattice of Hodge tensors, the image of $\pi_{1}(S)$ lies in $\mathbf{G}_{S}$ after replacing $S$ with a finite \'etale covering. Choosing a lattice $\Gamma \subset \mathbf{G}_{S}(\mathbb{Q})$ containing the image of monodromy, one obtains a \emph{period map} $\varphi : S \to \Gamma \backslash D_{S}$. We write $\pi : D_{S} \to \Gamma \backslash D_{S}$ for the obvious quotient map.
\begin{defn}
A special subvariety $Z \subset S$ is said to be of positive period dimension if $\dim \varphi(Z) > 0$. In other words, $\varphi(Z)$ is not a point.
\end{defn}
Now let $\mathcal{A} \subset \mathcal{S}$ denote the subset of maximal atypical special subvarieties, and further partition $\mathcal{A}$ as $\mathcal{A} = \mathcal{A}_{0} \sqcup \mathcal{A}_{\textrm{pos}}$, where $\mathcal{A}_{0}$ is the subset of period dimension zero and $\mathcal{A}_{\textrm{pos}}$ is the subset of positive period dimension. As the union of all subvarieties in $\mathcal{A}_{\textrm{pos}}$ is contained in a finite collection of strict special subvarieties by \autoref{geoZP}, one might expect that the main difficulty in resolving \autoref{ZPconj} is to control the subset of $\mathcal{A}_{0}$ which does not lie in a strict special subvariety of $S$. We therefore focus in on a particular subset of $\mathcal{A}_{0}$ defined by
\[ \mathcal{A}^{\textrm{iso}}_{0} := \{ Z \in \mathcal{A}_{0} : Z \textrm{ is a maximal strict special subvariety of }S \} . \]
We call the varieties in $\mathcal{A}^{\textrm{iso}}_{0}$ the \emph{isolated} atypical special subvarieties of period dimension zero.

\subsection{Period Domains of Definition}

In the remainder of the paper we assume that $\mathbf{G}^{\textrm{ad}}_{S}$ is $\mathbb{Q}$-simple, sometimes without mentioning it. We moreover fix a representation $\rho : \mathbf{G}_{S} \to \GL(V)$ coming from $\mathbb{V}$, where $V = \mathbb{V}_{s_{0}}$ is a fibre of $\mathbb{V}$, the choice of which is unimportant. We may thereby consider $D = D_{S}$ to be a collection of Hodge structures on $V$, and up to the action of $\Gamma$ view each period domain $D_{Z}$ associated to $Z \in \mathcal{S}$ as a complex submanifold of $D$. 

\begin{defn}
A \emph{period subdomain} $D' \subset D$ is a complex submanifold of the form $M_{h}(\mathbb{R}) \cdot h \subset D$, where $h \in D$ is a Hodge structure with Mumford-Tate group $M_{h}$. 
\end{defn}

\begin{defn}
\label{definedbydef}
A special subvariety $Z$ is said to be \emph{defined by} a Mumford-Tate subdomain $D'$ if it is an irreducible component of the inverse image of $\varphi^{-1}(\pi(D'))$. It is additionally said to be atypical for $D'$ if $\codim_{S} Z < \codim_{D} D'$. We also call $D'$ a domain of definition for $Z$.
\end{defn}

\paragraph{Hyperelliptic Example:} We illustrate the above definition with an example. Consider the hyperelliptic curve
\begin{align*}
C_{t_{0}} : y^2 &= x^5 + 20 x^4 - 26 x^3 + 20x^2 + x \\
&= x(x - \alpha_{1})(x - \alpha_{2})(x - \alpha_{3})(x - \alpha_{4}) .
\end{align*}
According to \cite[pg. 8]{kedlaya2008hyperelliptic}, the Jacobian of $C_{t_{0}}$ splits up to isogeny as a product of two elliptic curves, one of which has complex multiplication. We can consider $C_{t_{0}}$ as a fibre of the family 
\begin{align*}
C_{t} : y^2 &= x(x - \alpha_{1})(x - \alpha_{2})(x - \alpha_{3})(x - t) 
\end{align*}
where $t$ ranges over the complement in $\mathbb{A}^1$ of some finite set $E$ where the discriminant of the right-hand side vanishes. We obtain a family $f : C \to T := \mathbb{A}^1 - E$ of hyperelliptic curves, and then a variation of integral Hodge structure $\mathbb{V} := R^1 f_{*} \mathbb{Z}$ which is naturally polarized by the cup product. The monodromy representation associated to $\mathbb{V}$ has image equal to a finite index subgroup of $\Sp_{4}(\mathbb{Z})$ by \cite{yu1997toward} (cf. \cite[\S5]{zbMATH05224877}). The corresponding period map is a map $\varphi : S \to \mathcal{A}_{2} = \Sp_{4}(\mathbb{Z}) \backslash \mathbb{H}_{2}$, with $(\mathbf{G}, D) = (\textrm{Sp}_{4}, \mathbb{H}_{2})$.

Let $\widetilde{t}_{0} \in \mathbb{H}_{2}$ be a point lifting $\varphi(t_{0})$. Then the Hodge structure $h_{0}$ corresponding to $\widetilde{t}_{0}$ decomposes over $\mathbb{Q}$ as a direct sum $h_{0,\mathbb{Q}} = h_{1} \oplus h_{1}$, and we have a corresponding decomposition $V_{\mathbb{Q}} = V_{1} \oplus V_{2}$. Since only one of the summands has CM, say it is $V_{1}$, then there exists unique Hodge-theoretic idempotents $e_{1}$ and $e_{2}$ such that $e_{i} : V_{\mathbb{Q}} \to V_{i}$ is a surjective map of Hodge structures and $1 = e_{1} + e_{2}$. The summand $V_{1}$ also carries an additional Hodge-theoretic endomorphism $\eta : V_{1} \to V_{1}$ such that $\eta^2 = -d$, with $d > 0$ a positive integer. 

One can now consider the period subdomains $D_{e_{1}, \eta} \subset D_{e_{1}} \subset \mathbb{H}_{2}$ passing through $\widetilde{t}_{0}$. The first is the one-dimensional period subdomain where both $e_{1}$ and $\eta$ remain Hodge endomorphisms, and the second is the two-dimensional period subdomain where $e_{1}$ remains a Hodge idempotent. They are orbits of the real points of $\mathbb{Q}$-groups isogenous to $\SL_2 \times SO_2$ and $\SL_2 \times \SL_2$, respectively. Both period subdomains define the point $t_{0} \in S$, but only the first one does so atypically, since
\[ \codim_{S} \{ t_{0} \} = 1, \hspace{2em} \codim_{\mathbb{H}_{2}} D_{e_{1}, \eta} = 2, \hspace{2em}  \codim_{\mathbb{H}_{2}} D_{e_{1}} = 1 . \]

\begin{defn}
A period subdomain $D' \subset D$ is said to be \emph{defined by} a tensor $t \in \bigoplus_{a \geq 0} V^{\otimes a}$ if the Mumford-Tate group of a very general Hodge structure $h \in D'$ is exactly the stabilizer of $t$.
\end{defn}

\subsection{Complexity of Atypical Points}

We now assume $S$ is quasi-projective and count the varieties in $\mathcal{A}^{\textrm{iso}}_{0}$. For each integer $r$, we write $\mathcal{A}^{\textrm{iso}}_{0}(r) \subset \mathcal{A}^{\textrm{iso}}_{0}$ for the subset of special subvarieties which are defined atypically by a period domain which is itself defined by a tensor $t$ in $\bigoplus_{a \leq r} V^{\otimes a}$. For a positive integer $q$, we additionally write $\mathcal{A}^{\textrm{iso}}_{0}(r,q)$ for the further subset satisfying the property that the tensor $t$ can be chosen such that $Q(t,t) \leq q$. Note that $\mathcal{A}^{\textrm{iso}}_{0}(r,q)$ is a finite set. We also use the analogous notation $\mathcal{A}_{0}(r)$ and $\mathcal{A}_{0}(r,q)$ to consider the same question but considering all maximal atypical special varieties of zero period dimension, not necessarily isolated. After fixing a projective compactification $\overline{S}$ with ample line bundle $\mathcal{L}$, we can assign to each subvariety $Z \subset S$ a number $\deg Z := \deg_{\mathcal{L}} \overline{Z}$, with $\overline{Z}$ the closure of $Z$ in $\overline{S}$. We then write $\# \mathcal{A}^{\textrm{iso}}_{0}(r,q)$ for the sum of the degrees of the elements of $\mathcal{A}^{\textrm{iso}}_{0}(r,q)$, and likewise with $\# \mathcal{A}_{0}(r,q)$. In what follows, we always fix a compactification $S \subset \overline{S}$ with ample line bundle $\mathcal{L}$, the choice of which will be unimportant.

\begin{thm}
\label{complexitythm}
Let $(\mathbb{V}, Q)$ be an integral polarized variation of Hodge structure on a smooth quasi-projective complex algebraic variety $S$. Assume $\mathbf{G}_{S}^{\textrm{ad}}$ is $\mathbb{Q}$-simple. Let $r \geq 0$ be an integer. Then
\begin{align*}
\# \mathcal{A}^{\textrm{iso}}_{0}(r,q) = O(q^{\ep})
\end{align*}
for any $\ep > 0$. 
\end{thm}

This result can be applied to a conjecture of Grimm and Monnee. We recall from \cite[Def. 4.15]{BKU} that the \emph{level} of the variation of Hodge structure $\mathbb{V}$ is the unique integer $k$ such that the Hodge structure on the Lie algebra $\mathfrak{g}_{S}$ of $\mathbf{G}_{S}$ has exactly $2k+1$ non-zero Hodge summands. Using this, Grimm and Monnee predict in \cite[Conjecture 1]{grimm2024finitenesstheoremscountingconjectures} that:

\begin{conj}[Grimm-Monnee]
\label{GMconj}
Let $(\mathbb{V}, Q)$ be a variation of polarized integral Hodge structures of even weight $2k$ on a smooth quasi-projective variety $S$, and suppose that the level of $\mathbb{V}$ is at least three. Then 
\[ \# \mathcal{A}_{0}(1,q) = O(q^{\ep}) \]
for any $\ep > 0$. 
\end{conj}
As a corollary of \autoref{complexitythm} we obtain:

\begin{cor}
\label{conjholdscor}
The conjecture \autoref{GMconj} of Grimm-Monnee holds when $\mathbf{G}^{\textrm{ad}}_{S}$ is $\mathbb{Q}$-simple.
\end{cor}

\noindent Note that the simplicity assumption on $\mathbf{G}^{\textrm{ad}}_{S}$ almost always holds in explicit examples. Since by the main result of \cite{Andre1992} the derived subgroup of $\mathbf{G}_{S}$ contains a finite index subgroup of the image of monodromy, this assumption is satisfied whenever one has a large monodromy theorem, as for instance in the case of universal families of hypersurfaces and complete intersections in projective space \cite{zbMATH05651243}. 

One may wonder why one would conjecture \autoref{GMconj}, given that the Zilber-Pink conjecture \autoref{ZPconj} together with the atypicality of Hodge loci of level $3$ proven in \cite[Thm 3.3]{BKU} suggests that in fact $\mathcal{A}_{0}$ should be finite. It was explained to the author by the first author of \cite{grimm2024finitenesstheoremscountingconjectures} that he expects in many situations that the function $\# \mathcal{A}_{0}(1,q)$ of $q$ has a natural ``small $q$'' regime of initial growth where it behaves like an increasing subpolynomial function, before tapering off and becoming constant for large $q$. As a sanity check that such an expectation is not unreasonable one can then ask whether it is possible to show that the growth of $\# \mathcal{A}_{0}(1,q)$ is at most subpolynomial, which is what we do in this paper.

\subsection{Acknowledgements}

The author thanks Thomas Grimm for encouraging him to write up a proof of \autoref{complexitythm}, as well as Gregorio Baldi and Gal Binyamini for helpful discussions. 

\section{Recollections}

We collect various results from Hodge and unlikely intersection theory. We start with \cite[Thm 1.5]{defpermap}
\begin{thm}[BKT]
Given a period map $\varphi : S \to \Gamma \backslash D$, there exists a finite definable cover $S = \bigcup_{i=1}^{n} B_{i}$ by definable simply-connected open subsets and $\mathbb{R}_{\textrm{an,exp}}$-definable local lifts
\begin{align*}
\psi_{i} : B_{i} \to \mathfrak{O}_{i} \subset D
\end{align*}
of $\varphi$, where the $\mathfrak{O}_{i}$ are Hodge-theoretic Siegel sets. 
\end{thm}
We refer to \cite{defpermap} and \cite[\S2.1]{2024arXiv241208924U} for an introduction to Hodge-theoretic Siegel sets. For the purposes of reading this paper, the reader need not know the details of the definition. For a vector $v \in \mathbb{Q}^{m}$ for some $m$, its naive height is the maximum of the sizes of the numerators and denominators appearing the entries of $v$. If $V$ is integral, this is just the magnitude of the largest entry. The crucial property we will need is the following, which is a consequence of our work in \cite{2024arXiv241208924U}.
\begin{lem}
\label{heightsboundedbyQ}
Let $D$ be any period space for Hodge structures on a polarized lattice $(V, Q)$, and let $\mathfrak{O} \subset D$ be a Hodge-theoretic Siegel set. Then the naive heights of integral Hodge vectors $v \in V$ associated to points $x \in \mathfrak{O}$ are bounded by a polynomial in $Q(v,v)$. 
\end{lem}

\begin{proof}
By replacing $V$ with $V \otimes V$ and the vectors $v$ with $v \otimes v$, and restricting to just those Hodge vectors which are pure tensors, we may assume that $Q(v,v)$ is always a perfect square (cf. the proof in \cite[\S5.4]{2024arXiv241208924U}). Then by making the replacement $v \mapsto v/\sqrt{Q(v,v)}$, we can reduce to considering rational Hodge vectors $v \in V_{\mathbb{Q}}$ where $Q(v,v) = 1$, and bounding the heights of those vectors by a polynomial in their denominators. This is now immediate from \cite[Proposition 3.3]{2024arXiv241208924U}. 
\end{proof}

We will also have need for ``compact duals'': 

\begin{defn}
The compact dual $\ch{D}$ of $D$ is the orbit under $\mathbf{G}(\mathbb{C})$ of any point $h \in D$ inside the variety of all flags on $V_{\mathbb{C}}$ consisting of subspaces with the same dimensions as those of the Hodge flags parameterized by $D$. 
\end{defn}

\begin{defn}
Given a period subdomain $D' = M_{h}(\mathbb{R}) \cdot h \subset D$, its compact dual is $\ch{D}' := M_{h}(\mathbb{C}) \cdot h$. This is the Zariski closure of $D'$ in $\ch{D}$, and $D'$ is open in $\ch{D}'$. 
\end{defn}

The following is proven in \cite[Lem 7.12]{zbMATH08109694}, based on \cite[Thm 4.14]{voisin2010hodge}:

\begin{prop}
\label{conjclassesofMTgroups}
Hodge subdata $(M, D_{M}) \subset (\mathbf{G}, D)$ belong to finitely many $\mathbf{G}(\mathbb{R})$-equivalence classes, where the action on the first entry is by conjugation and the action on the second is by translation.
\end{prop}

Recall that for a definable subset $A \subset \mathbb{R}^{n}$ we define $A^{\textrm{alg}} \subset A$, the algebraic part, to be the union of all real semi-algebraic curves contained in $A$, and set $A^{\textrm{tr}} = A \setminus A^{\textrm{alg}}$. For a definable set $A$, we write $A[d, T]$ to be the set of points contained in $A$, defined over a number field of degree at most $d$, and with na\"ive height at most $T$. As a consequence of the Northcott property, this set is finite, and we write $\# A[d,T]$ for its cardinality. We have the following important result of Pila-Wilkie \cite{zbMATH05043321}, in the strengthened form proven by Pila in \cite[Thm. 5.3]{zbMATH05680945}.

\begin{defn}
A block in $\mathbb{R}^{m}$ is the image under a semi-algebraic map $\phi : \mathbb{R}^{n} \to \mathbb{R}^{m}$ of a connected open definable submanifold of a regular semi-algebraic set in $\mathbb{R}^{n}$ (cf. \cite[Def 3.2, 1]{zbMATH05680945}).
\end{defn}

\begin{thm}
\label{PW}
For any definable set $A \subset \mathbb{R}^{n}$, integer $d$, and constant $\ep > 0$, there exists a positve integer $J = J(A, d, \ep)$ and definable families of blocks $A_{j} \subset \mathbb{R}^{n} \times \mathbb{R}^{\mu_{j}}$, $j = 1, \hdots, J$, such that:
\begin{itemize}
\item[(1)] For each $j$ and $\eta \in \mathbb{R}^{\mu_{j}}$, we have $A_{j,\eta} \subset A$.
\item[(2)] The set $A[d,T]$ is contained in $O_{A,d,\ep}(T^{\ep})$ blocks contained in $A$, each a fibre of one of the families $A_{j}$; in particular,
\[ A^{\textrm{tr}}[d, T] = O_{A,d,\ep}(T^{\ep}) . \]
\item[(3)] Let $W \subset \mathbb{R}^{n}$ be the union over all $j$ of all the fibres of the $A_{j}$ of positive dimension. Then $W$ is definable, contained in $A^{\textrm{alg}}$, and
\[ (A - W)[d, T] = O_{A,d,\ep}(T^{\ep}) . \]
\end{itemize}
\end{thm}

\section{Proofs}

\subsection{Proof of \autoref{complexitythm}}

Using the main result of \cite{Andre1992}, the algebraic monodromy group $\mathbf{H}_{S}$ of $S$ (cf. \cite[Def. 1.15]{2024arXiv241208924U}) is a normal subgroup of the derived subgroup of $\mathbf{G}_{S}$. The fact that we assume $\mathbf{G}^{\textrm{ad}}_{S}$ is $\mathbb{Q}$-simple then means that $\mathbf{H}_{S}$ is either trivial or equal to the derived subgroup of $\mathbf{G}_{S}$. In the former case the period map $\varphi$ is constant (cf. \cite[III.A.2]{GGK}) and there is nothing to show, so we assume the latter. Then $\ch{D}$ is a homogeneous space for $\mathbf{H}_{S}$. 

One easily reduces to considering just those elements of $\mathcal{A}^{\textrm{iso}}_{0}(r)$ with dimension $d$, and then by intersecting $S$ with generic hyperplane sections and arguing as in \cite[\S5.1]{2024arXiv241208924U} we may reduce to the case where $d = 0$. We therefore write $\mathcal{P}^{\textrm{iso}}_{0} \subset \mathcal{A}^{\textrm{iso}}_{0}$ for the subset of maximal atypical varieties of zero dimension (instead of just zero period dimension), and analogously for $\mathcal{P}^{\textrm{iso}}_{0}(r)$ and $\mathcal{P}^{\textrm{iso}}_{0}(r,q)$. Note we can now identify $\mathcal{P}^{\textrm{iso}}_{0}$ and all of its subsets with subsets of $S(\mathbb{C})$, and we will do so; the quantities $\# \mathcal{P}^{\textrm{iso}}_{0}(r,q)$ are then just the cardinalities of these sets. Considering the definable cover $S = \bigcup_{i=1}^{n} B_{i}$, it suffices to reduce to proving the same bound for the subsets $\mathcal{P}_{i}(q) \subset \mathcal{P}^{\textrm{iso}}_{0}(r,q)$ corresponding to those points that lie in $B_{i}$, and then by computing the image $\psi_{i}(B_{i})$ and using that the number of isolated points in the fibres of the definable map $\psi_{i}$ is universally bounded, to proving the same bound for the subsets $\mathcal{B}_{i}(q) := \psi_{i}(\mathcal{P}_{i}(q)) \subset \psi_{i}(B_{i})$. We likewise set $\mathcal{B}_{i} := \psi_{i}(\mathcal{P}_{i})$. After replacing $\mathbb{V}$ with $\mathbb{V} \otimes \mathbb{V}$ if necessary we may assume that the Hodge structures parameterized by $\mathbb{V}$ have even weight. Twisting if necessary, we may then assume that the Hodge structures parameterized by $\mathbb{V}$ have weight zero. 

\begin{defn}
We say that a point $h \in \mathcal{B}_{i}$ is \emph{defined by} a Hodge subdatum $(M, D') \subset (\mathbf{G}, D)$ if $(M_{h}, M_{h}(\mathbb{R}) \cdot h) \subset (M, D')$ and $h$ is an isolated point of $\psi(B_{i}) \cap D'$. 
\end{defn}

Via \autoref{conjclassesofMTgroups} above, we have finitely many $\mathbf{G}_{S}(\mathbb{R})$-equivalence classes of Hodge subdata $(M, D_{M}) \subset (\mathbf{G}, D)$. Call such an equivalence class a \emph{type}, and label it by $\tau$. Because there are finitely many types, it suffices to prove the statement just for those elements of $\mathcal{B}_{i}$ for which the period domain which defines it atypically has type $\tau$. More precisely, we may consider
\begin{align*}
\mathcal{B}_{i}(\tau) :=& \left\{ h \in \mathcal{B}_{i} : \substack{ \textrm{exists }(M, D') \subset (\mathbf{G}, D)\textrm{ of type }\tau \textrm{ defining }h\textrm{ atypically} \\ \textrm{ and such that }D'\textrm{ is defined by some }t \in \bigoplus_{a \leq r} V^{\otimes a} } \right\} \\
\mathcal{B}_{i}(\tau, q) :=& \left\{ h \in \mathcal{B}_{i} : \substack{ \textrm{exists }(M, D') \subset (\mathbf{G}, D)\textrm{ of type }\tau \textrm{ defining }h\textrm{ atypically and such } \\ \textrm{ that }D'\textrm{ is defined by some }t \in \bigoplus_{a \leq r} V^{\otimes a} \textrm{ with }Q(t,t) \leq q } \right\}
\end{align*}
We may then reduce to showing that $|\mathcal{B}_{i}(\tau, q)| = O(q^{\ep})$ for any $\ep > 0$. We therefore drop the subscript and just write $\mathcal{B}(\tau) = \mathcal{B}_{i}(\tau)$, $\mathcal{B}(\tau, q) = \mathcal{B}_{i}(\tau, q)$ and $\psi(B) = \psi_{i}(B_{i})$. 

We write $\ch{D}$ for the compact dual of $D$; one can view $\ch{D}$ as the $G(\mathbb{C})$-orbit of any point $h \in D$ inside the algebraic variety of all flags on $V_{\mathbb{C}}$ consisting of subspaces with the same dimensions as those of the Hodge flags parameterized by $D$. Note that $D$ is an open subset of $\ch{D}$. Write $W_{a}$ for the affine space associated to the lattice $V^{\otimes a}$. The naive height on $V$ induces a naive height on each $V^{\otimes a}$, and these heights can be extended to a Weil heights $\theta_{a} : W_{a}(\overline{\mathbb{Q}}) \to \mathbb{R}_{\geq 0}$ and more generally a Weil height $\theta : \prod_{a \leq r} W_{a}(\overline{\mathbb{Q}}) \to \mathbb{R}_{\geq 0}$. One likewise obtains, using the polarization-induced isomorphism $V \simeq V^{\vee}$, heights on the affine spaces associated to $V^{\vee}$ and the tensor powers $V^{\otimes a} \otimes (V^{\vee})^{\otimes b}$, and then a height on $\mathbf{G}_{S}$ which one regards inside the affine space associated to $V \otimes V^{\vee}$. 

In what follows we fix, for each $\xi \in \mathcal{B}(\tau)$, the data of:
\begin{itemize}
\item[-] a pair $(M_{\xi}, D_{\xi})$ of type $\tau$ which defines $\xi$ atypically; and
\item[-] a tensor $t_{\xi} \in \bigoplus_{a \leq r} V^{\otimes a}$ defining $D_{\xi}$
\end{itemize}
subject to the condition that $Q(t_{\xi}, t_{\xi})$ is minimized. 

\begin{notn}
For a tensor $t$ we write $\ch{\NL}_{t} \subset \ch{D}$ for the locus of flags $F^{\bullet}$ such that $t \in F^{0}$, and $\mathbf{M}_{t} \subset \mathbf{G}$ for the subgroup stabilizing $t$. Given a Hodge datum $(M, D')$, we write $\ch{\NL}_{M} \subset \ch{D}$ for the locus of all flags $F^{\bullet}$ such that all tensors stabilized by $M$ lie in $F^{0}$. 
\end{notn}

Fix a representative $(M, D')$ for the type $\tau$, and let $\mathbf{L} \subset \mathbf{G}$ be the subgroup which stabilizes $\ch{D}'$ as a subvariety (each point of $\ch{D}'$ is sent to another point of $\ch{D}'$). Note that $M \subset \mathbf{L}$. 

\begin{prop}
\label{gconstrprop}
Fix a quasi-projective embedding $\iota : \mathbf{G} / \mathbf{L} \hookrightarrow \mathbb{P}^r$, and consider the induced Weil height $\theta_{\mathbf{L}} : (\mathbf{G}/\mathbf{L})(\overline{\mathbb{Q}}) \to \mathbb{R}_{\geq 0}$. Then there exists a map $\beta : \mathcal{B}(\tau) \to (\mathbf{G}/\mathbf{L})(\overline{\mathbb{Q}})$ with the following properties:
\begin{itemize}
\item[(i)] there exists a constant $d$ such that each point in the image of $\beta$ lies inside a number field of degree at most $d$;
\item[(ii)] for each $\xi \in \mathcal{B}(\tau)$ the height $\theta_{\mathbf{L}}(\beta(\xi))$ is bounded by a polynomial in $Q(t_{\xi}, t_{\xi})$; and
\item[(iii)] for each $\xi \in \mathcal{B}(\tau)$ we have $\beta(\xi) \cdot \ch{D}' = \ch{D}_{\xi}$.
\end{itemize}
\end{prop}

\begin{proof}
The proof is very similar to \cite[Prop. 7.15]{zbMATH08109694}. For illustrative purposes we divide into two cases. In what follows we write $W = \prod_{a \leq r} W_{a}$, which is the affine space corresponding to $\bigoplus_{a \leq r} V^{\otimes a}$.

\paragraph{$\dim D' = 0$:} In this case $\mathbf{G}/\mathbf{L} = \ch{D}$, and the points in $\mathcal{B}(\tau)$ all correspond to Hodge structures with complex multiplication. We consider the algebraic family $y : \mathcal{T} \to W$ whose fibres are the loci $\NL_{t}$. Then there is a constructible sublocus $W_{0} \subset W$, defined over $\mathbb{Q}$, such that the fibres of $y$ above $W_{0}$ are zero dimensional. We then consider the family $y_{0} : \mathcal{T}_{0} \to W_{0}$, where $\mathcal{T}_{0} = y^{-1}(W_{0})$. Then $y_{0}$ is quasi-finite, so the heights of the points in the fibres are polynomially bounded by the heights of the points in the image. In particular for each $t \in W_{0}$ associated to a point of $\mathcal{B}(\tau)$ we obtain that the associated point in $y^{-1}_{0}(t_{\xi})$ corresponding to $D_{\xi}$ has uniformly bounded field of definition and height bounded by a polynomial in the height of $t_{\xi}$, and hence a polynomial in $Q(t_{\xi},t_{\xi})$ using \autoref{heightsboundedbyQ}. Using the identification $\mathbf{G}/\mathbf{L} = \ch{D}$, properties (i), (ii), and (iii) follow immediately. 

\paragraph{$\dim D' > 0$:} Arguing in a similar fashion to \cite[Prop. 7.15]{zbMATH08109694}, there exists a finite union of locally closed $\mathbb{Q}$-subvarieties $\mathcal{T}_{1}, \hdots, \mathcal{T}_{\ell} \subset W \times \GL_{m}$ such that the union $\bigcup_{i} \mathcal{T}_{i}$ is exactly the locus of $(t, g)$ satisfying property that $\ch{\NL}_{t} = g\, \ch{\NL}_{M}$. Indeed, this is a constructible algebraic condition over $\mathbb{Q}$. Note that in the situation where $t$ is a Hodge tensor for some Hodge structure $h \in D$ with $M_{h} = \mathbf{M}_{t}$, and $(M_{h}, M_{h}(\mathbb{R}) \cdot h)$ has type $\tau$, then one can take for $g$ any $g \in \mathbf{G}(\mathbb{R})$ such that $g (M, D') = (M_{h}, M_{h}(\mathbb{R}) \cdot h)$ and the corresponding point $(t, g)$ lies inside some $\mathcal{T}_{i}$. 

By considering the projections to $W$, we obtain a finite collection $y_{i} : \mathcal{T}_{i} \to \mathcal{W}_{i}$ of algebraic families, with $\mathcal{W}_{i} \subset W$, such that the fibre of $y_{i}$ above $t$ consists of those $g$ satisfying (a) and (b). Each fibre of $y_{i}$ is naturally a torsor under the algebraic group $\mathbf{N}$ which preserves the algebraic variety $\ch{\NL}_{M}$, which by an analogous argument to \cite[VI.A.3]{GGK} (where $\mathbf{G} = \textrm{Aut}(V, Q)$) is just the normalizer of $M$ in $\mathbf{G}$. Fix some $y = y_{i}$ which we write as $y : \mathcal{T} \to \mathcal{W}$.
\begin{quote}
\textbf{Claim: }It suffices to show that, for any $t \in \mathcal{W}(\mathbb{Q})$, we can construct some $g_{t} \in y^{-1}(t)$ defined over a number field of uniformly bounded degree whose height is bounded by a uniform polynomial in in the height $\theta(t)$ of $t$.
\end{quote}
Indeed, let us suppose can achieve this, and fix representatives $n_{1}, \hdots, n_{k}$ defined over a number field for the components of $\mathbf{N} / M$. Considering some $\xi \in \mathcal{B}(\tau)$ we may set $\beta(\xi) := g_{t_{\xi}} \mathbf{L}$ and write $g_{\xi} = g_{t_{\xi}}$. Then (i) is true by assumption. We claim that (iii) is achieved up to replacing $g_{\xi}$ with $g_{\xi} n_{i}$ for some $n_{i}$ as above: indeed, by our reasoning above any $g \in \mathbf{G}(\mathbb{R})$ such that $g (M, D') = (M_{\xi}, M_{\xi}(\mathbb{R}) \cdot \xi)$ will produce a point in $y^{-1}(t_{\xi})$ which satisfies (iii). On the other hand action of $\mathbf{N}$ on the set of components of $\ch{\NL}_{M}$ factors through $\mathbf{N}/\mathbf{N}^{\circ}$, where $\mathbf{N}^{\circ} \supset M$ is the identity component of $\mathbf{N}$, so it suffices to choose $n_{i}$ such that $g_{\xi} n_{i}$ and $g$ lie in the same component of the torsor $y^{-1}(t_{\xi})$. Finally, (ii) is automatic from the properties of heights under polynomial maps together with \autoref{heightsboundedbyQ}. 

To show the claim, observe that $\mathbf{N}$ acts freely transitively on the fibres of $y$. Thus, after further stratifying $\mathcal{W}$ (cf. our argument in \cite[Prop. 7.15]{zbMATH08109694}), the map $y$ is an fppf torsor, and hence an \'etale torsor by \cite{localtorsorMO}. This implies that we can, after replacing $\mathcal{W}$ with an \'etale cover, choose an algebraic section of $y$. The result follows.
\end{proof}
For a definable set $A$, write $A^{0}$ for its subset of zero-dimensional components (isolated points). We consider the definable set
\begin{align*}
\mathcal{G} := \{ (F^{\bullet}, g) \in \ch{D} \times (\mathbf{G}/\mathbf{L})(\mathbb{C}) : F^{\bullet} \in [g \ch{D}' \cap \psi(B)]^{0} \} .
\end{align*}

By construction, the sets $\mathcal{G} \cap [\ch{D} \times \{ g \}]$ are finite for each $g \in (\mathbf{G}/\mathbf{L})(\mathbb{C})$, so applying \cite[Lem. 4.2]{2024arXiv241208924U} there exists a definable partition $\mathcal{G} = \mathcal{Y}_{1} \sqcup \cdots \sqcup \mathcal{Y}_{k}$ with the property that $\rho_{i} : \mathcal{Y}_{i} \to (\mathbf{G}/\mathbf{L})(\mathbb{C})$ is injective for each $i$. If we write $\mathcal{I}_{i}$ for the image of $\rho_{i}$ and $\mathcal{D}_{i} \subset \ch{D}$ for the projection of $\mathcal{Y}_{i}$ to $\ch{D}$, then we obtain definable maps $e_{i} : \mathcal{I}_{i} \to \psi(B) \subset \mathcal{D}_{i}$ which are obtained as compositions $\textrm{pr}_{\ch{D}} \circ \rho^{-1}_{i}$. Note that, for each $\xi \in \mathcal{B}(\tau)$, the point $(\xi, \beta(\xi))$ is a point of some $\mathcal{Y}_{i}$. 

\vspace{0.5em}

Now suppose that, for some constant $c > 0$, we have that $|\mathcal{B}(\tau,q)| \geq c q^{\ep}$ for infinitely many $q$ as $q \to \infty$. Then necessarily, after possibly shrinking $c$, there is some $i$ such that
\begin{equation}
\label{newsettobound}
| \underbrace{\{ (\xi, \beta(\xi)) : \xi \in \mathcal{B}(\tau,q) \} \cap \mathcal{Y}_{i}}_{\mathcal{C}_{i}(q)} | \geq c q^{\ep}
\end{equation}
for infinitely many $q$ as $q \to \infty$. We fix this index $i$ and write $\mathcal{Y} = \mathcal{Y}_{i}$, $\mathcal{I} = \mathcal{I}_{i}$, $\mathcal{C}(q) = \mathcal{C}_{i}(q)$ and $\mathcal{D} = \mathcal{D}_{i}$. Note that if we start with a point $(\xi, \beta(\xi)) \in \mathcal{C}(q)$ for $\xi \in \mathcal{B}(\tau)$, then $e(\beta(\xi)) = \xi$. 

Now fix, using \autoref{gconstrprop}(ii), an exponent $m$ and constant $\kappa$ such that $\theta_{\mathbf{L}}(\beta(\xi)) \leq \kappa\, Q(\xi, \xi)^{m}$ for all $\xi \in \mathcal{B}(\tau)$. Each point $\xi \in \mathcal{B}(\tau)$ such that $(\xi, \beta(\xi)) \in \mathcal{C}(q)$ satisfies the property that $\beta(\xi) \in \mathcal{I}[d, \kappa q^{m}]$ with $d$ as in \autoref{gconstrprop}(i), so necessarily as $\rho$ is injective
\begin{equation}
\label{newsettobound2}
|\mathcal{I}[d, \kappa q^{m}]| \geq c q^{\ep} 
\end{equation}
for infinitely many $q$ as $q \to \infty$, where we use that $\rho$ is injective. Then using \autoref{PW}(ii) and choosing the $\ep$ in \autoref{PW} to be our $\ep/2m$, the set $\mathcal{I}[d, \kappa q^{m}]$ is contained in $O(q^{\ep/2})$ blocks contained in $\mathcal{I}$. In particular, we can take $q$ large enough such that one such block, call it $A$, contains at least two points $\beta(\xi)$ and $\beta(\xi')$ with $\beta(\xi) \neq \beta(\xi')$. Then the definable function $e$ takes both the values $\xi$ and $\xi'$ on $A$, hence is non-constant on $A$. By the definition of block $A$ is connected, so the image $e(A)$ is a connected definable set containing at least two distinct points $\xi$ and $\xi'$, so necessarily of positive definable dimension.

By replacing $A$ with a semi-algebraic curve $A' \subset A$ which passes through $\beta(\xi)$ and $\beta(\xi')$, we may assume that $A$ has real dimension $1$. We can then choose some real semi-algebraic curve $\widetilde{A} \subset \mathbf{G}(\mathbb{C})$ which maps surjectively onto $A$ and thereby learn that the union 
\[ \bigcup_{a \in \widetilde{A}} a \ch{D}' \cap \psi(B) = (\widetilde{A} \cdot \ch{D}') \cap \psi(B) \]
has a component of positive definable dimension (i.e., the one containing $e(A)$) passing through both $\xi$ and $\xi'$. Applying \cite[Lem. 4.1]{zbMATH06144656} (cf. \cite[Lem. 7.18]{zbMATH08109694}) and the fact that $\mathcal{I}$ is complex-analytic, we may even assume that $\widetilde{A}$ is a complex-algebraic curve. Set $\ch{E} := \widetilde{A} \cdot \ch{D}'$. Let $C \subset \ch{E} \cap \psi(B)$ be a positive-dimensional complex analytic component containing $\xi$.

Now because $\widetilde{A}$ is a curve, the constructible algebraic set $\ch{E}$ has (complex) dimension $\dim \ch{D}' + 1$. Since the points in $\mathcal{B}(\tau)$ are defined atypically by translates of $\ch{D}'$, this means in particular that
\[ \codim_{\psi(B)} C = \dim \psi(B) - 1 < \codim_{\ch{D}} \ch{D}' - 1 = \codim_{\ch{D}} \ch{E} . \]
Applying the Ax-Schanuel Theorem for variations of Hodge structures \cite[Thm. 1.1]{AXSCHAN} together with the fact that $\ch{D}$ is equal to the $\mathbf{H}$-orbit of a point in $\ch{D}$ (note in \cite{AXSCHAN} that they write $\mathbf{G}$ for our $\mathbf{H}$), one learns that $\varphi^{-1}(\pi(C))$ lies in a finite union of strict positive-dimensional (in the sense of period dimension) weakly special subvarieties of $S$. By construction the special point in $\mathcal{A}^{\textrm{iso}}_{0}(r)$ mapping to $\xi$ lies inside one of these subvarieties. Moreover these weakly special subvarieties can be assumed to be special subvarieties of $S$ as a consequence of the $\mathbb{Q}$-simplicity of $\mathbf{G}^{\textrm{ad}} = \mathbf{H}^{\textrm{ad}}$ as well as \cite[Lem 2.5]{fieldsofdef}. But points of $\mathcal{A}^{\textrm{iso}}_{0}(r)$ do not lie in strict positive-dimensional special subvarieties of $S$ by assumption.

\vspace{0.5em}

We conclude that our assumption that $\mathcal{B}(\tau,q) \geq c q^{\ep}$ for some $c, \ep > 0$ as $q \to \infty$ is violated, so $\mathcal{B}(\tau, q) = O(q^{\ep})$ for any $\ep > 0$, completing the proof.

\begin{rem}
The second last paragraph of the above argument is the only portion that uses that $\mathbf{G}_{S}^{\textrm{ad}}$ is $\mathbb{Q}$-simple, which is used to ensure the notion of atypicality in terms of special varieties agrees with that coming from Ax-Schanuel, as well as to ensure that a strict positive-dimensional weakly special subvariety lies in a strict special subvariety.
\end{rem}

\subsection{Application to \autoref{GMconj}}

The statement \autoref{complexitythm} immediately specializes to \autoref{GMconj} by substituting $r = 1$, except for the possible difference between $\mathcal{A}^{\textrm{iso}}_{0}(1,q)$ and $\mathcal{A}_{0}(1,q)$. However we observe that:

\begin{lem}
When the level of $\mathbb{V}$ is at least three, then $\mathcal{A}_{0} = \mathcal{A}^{\textrm{iso}}_{0}$.
\end{lem}

\begin{proof}
If some variety $Z$ is in $\mathcal{A}_{0}$ but not in $\mathcal{A}^{\textrm{iso}}_{0}$ then it is a maximal atypical special subvariety of zero period dimension contained in some maximal strict special subvariety $Y$ of $S$ of positive period dimension (a priori not necessarily atypical). Since $\mathbb{V}$ has level at least three, all such $Y$ are in fact atypical by \cite[Thm 3.3]{BKU}, so this cannot happen because $Z \in \mathcal{A}_{0}$ means it is a maximal atypical special subvariety. 
\end{proof}

Note that in particular one obtains the following strengthening of \autoref{conjholdscor}:

\begin{thm}
\label{complexitythm2}
Let $(\mathbb{V}, Q)$ be an integral polarized variation of Hodge structure on a smooth quasi-projective complex algebraic variety $S$. Assume $\mathbf{G}_{S}^{\textrm{ad}}$ is $\mathbb{Q}$-simple and $\mathbb{V}$ has level at least three. Let $r \geq 0$ be an integer. Then
\begin{align*}
\# \mathcal{A}_{0}(r,q) = O(q^{\ep})
\end{align*}
for any $\ep > 0$. 
\end{thm}

\bibliography{hodge_theory}
\bibliographystyle{alpha}

\end{document}